\documentclass[12pt, draft]{amsart}
\usepackage{eucal}

\theoremstyle{plain}
\newtheorem*{main}{\textsf{Main Theorem}}
\newtheorem{theorem}{\textsf{Theorem}}
\newtheorem{lemma}{\textsf{Lemma}}
\newtheorem{corollary}{\textsf{Corollary}}

\theoremstyle{remark}
\newtheorem{remark}{Remark}

\newcommand{\A}{\mathcal{A}}
\newcommand{\BH}{\mathcal B(\mathcal H)}

\newcommand{\cH}{\mathcal H}
\newcommand{\CH}{\mathcal{C}(\mathcal{H})}
\newcommand{\T}{\mathfrak{T}}
\newcommand{\li}{\ell_\infty(\mathbb{N}, \mathcal B(\mathcal H))}
\newcommand{\N}{\mathbb{N}}

\begin{document}
\title[Reflexivity of the automorphism and isometry groups]{The
                    automorphism and isometry groups of\\
                    $\li$ are topologically reflexive}
\author{Lajos Moln\'ar}
\address{Institute of Mathematics\\
         Lajos Kossuth University\\
         4010 Debrecen, P.O.Box 12, Hungary}
\email{molnarl@math.klte.hu}
\thanks{This paper was written when the author,
        holding a scholarship of the Volkswagen-Stiftung
        of the Konferenz der Deutschen Akademien der Wissenschaften,
        was a visitor at the University of Paderborn, Germany.
        He is very grateful to
        Prof. K.-H.~Indlekofer for his kind hospitality.
        This research was partially supported also by the Hungarian
        National
        Foundation for Scientific Research (OTKA), Grant No. T--016846
        F--019322}
\keywords{Reflexivity, automorphism group, isometry group, operator
        algebra, Stone-\v Cech compactification}
\subjclass{Primary: 47B49, 46L40, 47D25; Secondary: 54D35}
\date{June 11, 1997}

\begin{abstract}
The aim of this paper is to show that the automorphism and isometry
groups of the $C^*$-algebra $\li$ are topologically reflexive which,
as one of our former results shows, is not the case with the "scalar
algebra" $\ell_\infty$.
\end{abstract}
\maketitle

\section{Introduction and Statement of the Results}

Let $\A$ be a $C^*$-algebra. A subset $\T$ of the algebra
$\mathcal{B}(\A)$ of all bounded linear transformations
on $\A$ is called topologically
reflexive if for every $\Phi\in \mathcal{B}(\A)$, the condition that
$\Phi(A)\in \overline{\T(A)}$ (the norm-closure of $\T(A)$)
holds true for every $A\in \A$ implies $\Phi\in \T$.
Similarly, we say that $\T$ is algebraically reflexive if we have
$\Phi\in \T$ for any $\Phi \in \mathcal{B}(\A)$ with $\Phi(A)\in \T(A)$
$(A\in \A)$. Roughly speaking, reflexivity means that the elements of
$\T$ are, in some sense, completely determined by their local actions.

In what follows, let $\cH$ stand for a complex separable infinite
dimensional Hilbert space. Reflexivity problems concerning subspaces of
$\BH$ represent one of the most active and important research areas in
operator theory. The question of reflexivity in our sense
described
above was first considered by Kadison \cite{Kad} and Larson and Sourour
\cite{LaSo} in the case of derivation algebras of operator algebras.
In \cite{BrSe} Bre\v sar and \v Semrl proved the algebraic reflexivity
of the automorphism group of $\BH$.
We emphasize that by an automorphism we mean a multiplicative linear
bijection, so we do not assume the *-preserving property.
As for topological reflexivity,
the first result was obtained by Shul'man in \cite{Shu} concerning the
derivation algebra of any $C^*$-algebra. In our papers \cite{BaMo,
MolStud1,
MolStud2} we investigated the topological reflexivity of the groups of
all *-automorphisms, respectively automorphisms, respectively surjective
isometries of some operator algebras. The results which are in the
closest relation to our present investigation are the following. In
\cite{MolStud1} we showed that the automorphism and isometry groups of
$\BH$
are topologically reflexive. In \cite{BaMo} we obtained that in the case
of the commutative algebra $\ell_\infty$, these groups are topologically
nonreflexive. The main result of this paper is
that the "mixed" $C^*$-algebra
\[
\li=\{ (A_n)\, :\, A_n\in \BH\,\, (n\in \N),\, \sup_n \| A_n\| <\infty
\}
\]
has topologically reflexive automorphism and isometry groups.
This could be surprising, since these groups in question "cannot be more
reflexive"
than they are in the case of $\mathbb{C}$, the role of which is played
by the operator algebra $\BH$ in $\li$. On the way to the
proof of this assertion we present some hopfully interesting
results concerning the tensor product of $\ell_\infty$ and $\BH$ as well
as the Stone-\v Cech compactification $\beta \N$ of $\N$.

In what follows we need the concept of Jordan homomorphisms.
If $\mathcal{R}$ and $\mathcal{R'}$ are algebras, then a linear map
$\phi :\mathcal{R} \to \mathcal{R'}$ is called a Jordan homomorphism if
\[
\phi(A)^2=\phi(A^2) \qquad (A\in \mathcal{R}).
\]
The following equations are well-known to be fulfilled by any Jordan
homomorphism
\begin{subequations}\label{E:jor}
\begin{gather}
\phi(A)\phi(B)+\phi(B)\phi(A)=\phi(AB+BA) \label{E:jor1}\\
\phi(A)\phi(B)\phi(A)=\phi(ABA)           \label{E:jor2}\\
\phi(A)\phi(B)\phi(C)+\phi(C)\phi(B)\phi(A)=
\phi(ABC+CBA)                             \label{E:jor3}
\end{gather}
\end{subequations}
where $A,B,C\in \mathcal{R}$ (e.g. \cite[6.3.2 Lemma]{Pal}).

It is well-known that the $C^*$-algebras $\ell_\infty$ and $C(\beta \N)$
(the algebra of all continuous complex valued functions on
$\beta \N$) are isomorphic. In fact, this follows
from the property of the Stone-\v Cech compactification that every
continuous function from a completely regular space $X$ into
a compact Hausdorff space
$Y$ can be uniquely extended to a continuous function defined on
$\beta X$. The map which sends every element of $\ell_\infty$ to its
unique
extension to an element of $C(\beta \N)$ gives the desired isomorphism.
Due to the topological properties of $\beta \N$, we have "singular"
characters of the
commutative algebra $\ell_\infty$ by simply considering any one-point
evaluation functional on
$C(\beta \N)$ which corresponds to a point in $\beta \N \setminus \N$.
The word "singular" means here that this character annihilates the
cofinite sequences in $\ell_\infty$ (and hence it is not
$w^*$-continuous). In our first theorem we show that
this is not the case if we replace the set $\mathbb{C}$ of values by
the operator algebra $\BH$.

\begin{theorem}\label{T:linosin}
There is no nonzero Jordan homomorphism
$\Phi:\li \to \BH$ which vanishes on the set of all cofinite sequences.
\end{theorem}

This result has the corollary that the kernels of irreducible Jordan
homomorphisms of the above type correspond to elements of $\N$.
Comparing this with the case of the characters of $\ell_\infty$, one can
say that the operator algebra $\BH$ cleans up
the Stone-\v Cech compactification of $\N$.

\begin{corollary}\label{C:linosin}
Let $\Phi :\li \to \BH$ be an irreducible Jordan homomorphism
(i.e. a Jordan homomorphism whose range has only trivial invariant
subspaces). Then there is a positive integer $n$ so that
\[
\ker \Phi =\{ (A_k)_k\in \li \, : \, A_n =0\}.
\]
\end{corollary}

It is well-known in the theory of tensor products that for any compact
Hausdorff space $X$ and for any $C^*$-algebra $\mathcal{A}$, the
$C^*$-algebras $C(X)\otimes \mathcal{A}$ and $C(X, \mathcal{A})$ (the
space of all continuous functions from $X$ into $\mathcal{A}$) are
isomorphic. The following corollary shows that a similar assertion does
not hold true if one considers $\ell_\infty$ instead of $C(X)$.

\begin{corollary}\label{C:linoiso1}
The $C^*$-algebras $\ell_\infty \otimes \BH$ and $\li$ are nonisomorphic
even as Jordan algebras.
\end{corollary}

This statement can be reformulated in the following way.

\begin{corollary}\label{C:linoiso2}
The $C^*$-algebras $\li$ and $C(\beta \N, \BH)$ are nonisomorphic
even as Jordan algebras.
\end{corollary}

One short remark should be added here. Namely, for any finite
dimensional
Hilbert space $\mathcal{K}$ we do have an isomorphism between the
$C^*$-algebras $\ell_\infty(\N, \mathcal{B}(\mathcal{K}))$ and $C(\beta
\N, \mathcal{B}(\mathcal{K}))$.
In fact, since in that case the functions in $\ell_\infty(\N,
\mathcal{B}(\mathcal{K}))$
have precompact ranges, one can use the function extension property of
the Stone-\v Cech compactification. Nevertheless,
Corollary~\ref{C:linoiso2} shows that not only this extension stuff
breaks down in the infinite dimensional case, but there does not exist
any isomorphism between $\li$ and $C(\beta \N, \BH)$.

In the next theorem which can also be consireded as a corollary of
Theorem~\ref{T:linosin}, we describe the automorphism and isometry
groups of $\li$.

\begin{theorem}\label{T:liaut}
Let $\Phi:\li \to \li$ be an automorphism. Then there are
automorphisms $\phi_n$ $(n\in \N)$ of $\BH$ and a bijection $\varphi:\N
\to \N$ so that $\Phi$ is of the form
\begin{equation}\label{E:liaut}
\Phi((A_k)_k)=(\phi_n(A_{\varphi(n)})) \qquad ((A_k)_k\in \li).
\end{equation}
Analogue statement holds true for the surjective isometries of $\li$
as well.
\end{theorem}

To be more specific with \eqref{E:liaut}, we recall
the well-known folk results that every automorphism of $\BH$ is of the
form \[
A \mapsto TAT^{-1}
\]
with some invertible operator $T\in \BH$ and that every surjective
isometry of $\BH$ is either of the form
\[
A\mapsto UAU^*
\]
or of the form
\[
A\mapsto UA^{tr}U^*
\]
with some unitary operator $U\in \BH$, where ${}^{tr}$ denotes the
transpose operation with respect to an arbitrary but fixed complete
orthonormal system in $\cH$.

After these preliminary results we shall be able to prove the main
result of the paper which follows.

\begin{main}\label{T:liref}
The automorphism and isometry groups of $\li$ are
topologically reflexive.
\end{main}

\section{Proofs}

We begin with a lemma which we shall use repeatedly throughout.

\begin{lemma}\label{L:jorker}
Let $\phi: \BH \to \BH$ be a Jordan homomorphism. Then $\phi$ is either
injective or we have $\phi=0$.
\end{lemma}

\begin{proof}
First observe that by \eqref{E:jor1} every Jordan homomorphism
preserves the idempotents as well as the orthogonality between
them (the idempotents $P,Q$ are called orthogonal if $PQ=QP=0$).
Now, since every Jordan homomorphism of $\BH$
is continuous (see \cite[Lemma 1]{MolStud1}), the kernel of $\phi$ is a
closed Jordan ideal of a $C^*$-algebra and hence it is an associative
ideal as well (see \cite{CiYo}). Therefore, $\ker \phi$ is either
$\{0\}$, $\CH$ (the ideal of all compact operators on $\cH$) or $\BH$.
Since the Calkin algebra $\BH/\mathcal{C}(\cH)$ has
uncountably many pairwise orthogonal nonzero idempotents which does not
hold true for $\BH$, we have the assertion.
\end{proof}

\begin{proof}[Proof of Theorem~\ref{T:linosin}]
Suppose that $\Phi:\li \to \BH$ is a Jordan
homomorphism which vanishes on the cofinite sequences.
Define $\phi: \BH \to \BH$ by
\[
\phi(A)=\Phi(A,A, \ldots ) \qquad (A\in \BH).
\]
Obviously, $\phi$ is a Jordan homomorphism.
Let $(e_n)$ be a complete orthonormal sequence in $\cH$ and denote
by $S\in \BH$ the corresponding unilateral shift.
For any $n\in \N$, let $P_n$ be the orthogonal projection onto the
subspace generated by $e_1, \ldots, e_n$. If $n,m\in \N$, since
$\Phi$ vanishes on the cofinite sequences, by \eqref{E:jor3} we
can compute
\begin{multline} \label{E:nosin1}
\phi(P_n)\Phi(S, S^2, S^3, \ldots )\phi(P_m)+
\phi(P_m)\Phi(S, S^2, S^3, \ldots )\phi(P_n)=\\
\Phi( P_nSP_m+P_mSP_n, P_nS^2P_m+P_mS^2P_n, P_nS^3P_m+P_mS^3P_n,
\ldots )=0.
\end{multline}
Since $\phi :\BH \to \BH$ is a Jordan homomorphism, it follows from
\cite[Lemma 2]{MolStud1} that the sequence $(\phi(P_n))$ converges
strongly to an idempotent $E\in \BH$ which does not depend on
the particular choice of $(e_n)$ and $E$ commutes
with the range of $\phi$. Then the map $\phi':\BH \to\BH$ defined by
$\phi'(A)=\phi(A)(I-E)$ is a Jordan homomorphism which vanishes on the
finite rank operators. By Lemma~\ref{L:jorker},
we have $\phi' =0$ and this gives us that
$\phi(.)=\phi(.)E=E\phi(.)$.
Now, using \eqref{E:jor2}, from \eqref{E:nosin1} we obtain
\[
0=\phi(I)E\Phi(S, S^2, S^3, \ldots )E\phi(I)=
\phi(I)\Phi(S, S^2, S^3, \ldots )\phi(I)=
\Phi(S, S^2, S^3, \ldots ).
\]
By \eqref{E:jor1} this further implies that
\begin{multline*}
0=
\Phi(S, S^2, S^3, \ldots )\Phi(S^*, {S^*}^2, {S^*}^3, \ldots )+
\Phi(S^*, {S^*}^2, {S^*}^3, \ldots )\Phi(S, S^2, S^3, \ldots )=\\
\Phi((I-P_1, I-P_2, I-P_3, \ldots )+(I,I,I, \ldots))
\end{multline*}
which results in
\begin{equation} \label{E:nosin2}
2\Phi(I,I,I, \ldots)=\Phi(P_1, P_2, P_3, \ldots ).
\end{equation}
Since the operators
$\Phi(I,I,I, \ldots)$ and $\Phi(P_1, P_2, P_3, \ldots )$ are
idempotents, we deduce from \eqref{E:nosin2} that $\Phi(I,I,I, \ldots
)=0$.
By \eqref{E:jor1} we clearly have $\Phi=0$. This completes the proof.
\end{proof}

\begin{remark}\label{R:linosin}
Observe that an inspection of the previous proof
(see \eqref{E:nosin1}) shows that the conclusion in
Theorem~\ref{T:linosin} remains true if we assume that
$\Phi:\li \to \BH$ vanishes only
on all cofinite sequences with finite-rank coordinates.
\end{remark}

\begin{proof}[Proof of Corollary \ref{C:linosin}]
First, by Theorem
\ref{T:linosin}, we obtain that $\Phi$ takes a nonzero value on a
cofinite sequence, say $(A_1, \ldots, A_n, 0, \ldots)$.
Therefore, in the equation
\[
0\neq \Phi(A_1, \ldots, A_n, 0, \ldots)=
\Phi(A_1, 0, \ldots)+ \dots +
\Phi(0, \ldots, 0, A_n, 0, \ldots)
\]
one term on the right hand side, say the last one, must be nonzero.
Thus, the map
\[
\Phi_n : A \mapsto \Phi (0, \ldots, 0, \underset{n}{A}, 0, \ldots)
\]
is a nonzero Jordan homomorphism of $\BH$. Now, by \eqref{E:jor1} it
is easy to see that the idempotent $\Phi_n(I)$ is necessarily nonzero.
Furthermore, using \eqref{E:jor2}, for any $(A_k)_k\in \li$ we infer
\[
2\Phi_n(I)\Phi((A_k)_k)\Phi_n(I)=2\Phi_n(A_n)=\Phi((A_k)_k)\Phi_n(I)+
\Phi_n(I)\Phi((A_k)_k).
\]
Multiplying this equation by $\Phi_n(I)$ from the right, we obtain that
\[
\Phi_n(I)\Phi((A_k)_k)\Phi_n(I)=\Phi((A_k)_k)\Phi_n(I)
\]
which means
that the range of $\Phi_n(I)$ is an invariant subspace of the range of
$\Phi$. Therefore, we have $\Phi_n(I)=I$. Now, we can compute
\[
2\Phi((A_k)_k)=
\Phi((A_k)_k)\Phi_n(I)+\Phi_n(I)\Phi((A_k)_k)=
2\Phi_n(A_n).
\]
Clearly, it remains to prove that $\Phi_n$ is injective. Since
$\Phi_n \neq 0$, by Lemma~\ref{L:jorker} we have $\ker \Phi_n=\{0\}$.
\end{proof}

\begin{proof}[Proof of Corollary~\ref{C:linoiso2}]
Let us suppose on the contrary that there exists a Jordan isomorphism
$\Phi: \li \to C(\beta \N, \BH)$.
Plainly, $\Phi$ preserves the nonzero minimal idempotents in these
algebras.
Any idempotent $f\in C(\beta \N, \BH)$ is a continuous function whose
values are idempotents
in $\BH$. Obviously, if we multiply $f$ by the characteristic
function of any point in $\N \subset \beta \N$, we obtain an
idempotent in $C(\beta \N, \BH)$. On the other hand, if a
function $f\in C(\beta \N, \BH)$ vanishes on $\N$, then by
$\| f(.)\| \in C(\beta \N)$ we have $f=0$. Putting these together,
we find that the nonzero minimal idempotents in $C(\beta \N, \BH)$ are
those
functions on $\beta \N$ which take only one nonzero value, they take it
at some point in $\N$ and this value is a rank-one idempotent in
$\BH$.
Now, let $p\in \beta \N \setminus \N$ and consider the homomorphism
$\Psi: C(\beta \N, \BH)\to \BH$ defined by $\Phi(f)=f(p)$.
Clearly, $f(p)=0$ for every minimal idempotent in $C(\beta \N, \BH)$.
This gives us that the nonzero Jordan homomorphism $\Psi \circ \Phi :\li
\to \BH$ vanishes on every cofinite sequence having finite rank
coordinates. But this is a contradiction by Remark~\ref{R:linosin}.
\end{proof}

\begin{proof}[Proof of Corollary \ref{C:linoiso1}]
Since $\ell_\infty$ and $C(\beta \N)$ are isomorphic, we obtain that
$\ell_\infty \otimes \BH$ and $C(\beta \N)\otimes
\BH \simeq C(\beta \N, \BH)$ are isomorphic as $C^*$-algebras.
Now, the statement follows from Corollary \ref{C:linoiso2}.
\end{proof}

\begin{proof}[Proof of Theorem \ref{T:liaut}]
First let $\Phi:\li \to \li$ be a Jordan automorphism.
Plainly, every coordinate function $\Phi^m$ of $\Phi$ satisfies
the condition in Corollary~\ref{C:linosin}. Hence, to every $m\in \N$
there corresponds a positive integer $\varphi(m)$ such that
\[
\Phi^m((A_k))=\Phi^m(0, \ldots, 0, A_{\varphi(m)}, 0, \ldots)=
\Phi^m_{\varphi(m)}(A_{\varphi(m)}),
\]
where $\Phi^m_{\varphi(m)}$ is a Jordan automorphism of $\BH$ (see the
proof of Corollary~\ref{C:linosin}). Let us show that $\varphi:\N \to
\N$ is a bijection.
Since $\Phi$ preserves the nonzero minimal idempotents in $\li$, thus
for any rank-one idempotent $P\in \BH$, two different
coordinates of $\Phi(0,\ldots, 0, P, 0, \ldots)$ cannot be nonzero.
This readily implies that $\varphi$ is injective.
Assume now that $\varphi$ is not surjective. Then it is easy to see that
there is at least one coordinate, say the first one, for which
$\Phi(A, 0, \ldots, )=0$ $(A\in \BH)$. But this contradicts
the injectivity of $\Phi$.
Consequently, we obtain that $\varphi$ is bijective.
Since we apparently have
\[
\Phi((A_k)_k)=(\Phi^m_{\varphi(m)}(A_{\varphi(m)})) \qquad ((A_k)_k\in
\li),
\]
the first statement of the theorem follows.
As for the second one, we refer to a well-known theorem of Kadison
stating that the surjective isometries of
a $C^*$-algebra are exactly those maps which can be written as a Jordan
*-automorphism multiplied by a fixed unitary element (see
\cite[7.6.16, 7.6.17]{KadRing2}).
\end{proof}

\begin{proof}[Proof of the Main Theorem]
Let first $\Phi :\li \to \li$ be a
continuous linear map which is an approximately local automorphism,
i.e. suppose that $\Phi$ can be approximated at every element
$(A_k)_k\in \li$ by the values of a sequence of automorphisms at
$(A_k)_k$.
Using Theorem~\ref{T:liaut}, one can readily verify that the map
\[
A\mapsto \Phi^m(A,A, \ldots )
\]
is an approximately local automorphism of $\BH$ and then, by
\cite[Theorem 2]{MolStud1}, it is an automorphism. In particular,
$\Phi^m$ is surjective and using Corollary~\ref{C:linosin} as well as
its proof, we obtain that for every $m\in \N$ there is
a positive integer $\varphi(m)$ so that $\Phi_{\varphi(m)}^m$ is an
automorphism and we have
\[
\Phi((A_k))=(\Phi^m((A_k)_k))=(\Phi^m_{\varphi(m)}(A_{\varphi(m)}))
\qquad ((A_k)_k\in \li).
\]
We show that the function $\varphi:\N \to \N$ is a bijection.
Like in the proof of Theorem~\ref{T:liaut}, we obtain immediately
that $\varphi$ is injective.
The surjectivity is also easy to see. Indeed, if
$\varphi$ were not surjective, then we would have
such a coordinate, say the first one for which
$\Phi(A, 0, \ldots)=0$ holds true for every $A\in \BH$.
But this is a contradiction, since, for example, the norm of
the image of $(I,0,0,\ldots)$ under any automorphism of $\li$ is 1.
For every $m\in \N$ let
$T_m\in \BH$ be an invertible operator such that
\[
\Phi^m_{\varphi(m)}(A)=T_mAT_m^{-1} \qquad (A\in \BH).
\]
Since $\Phi$ is of the form
\[
\Phi((A_k)_k)=(T_mA_{\varphi(m)}T_m^{-1}) \qquad ((A_k)_k\in \li)
\]
and $\Phi$ is bounded, considering sequences
of the form $(x_k \otimes y_k)_k$, one can easily verify
that $\sup_k \| T_k\|
\| T_k^{-1}\| <\infty$. This immediately gives us the bijectivity of
$\Phi$, concluding the proof of our statement in the case of
automorphisms.

Let next $\Phi$ be an approximately local surjective isometry of
$\li$. Since $\Phi$ clearly preserves the unitaries in the $C^*$-algebra
$\li$, it follows from a well-known theorem of Russo and Dye
\cite[Corollary 2]{RuDy} that $\Phi$ is a Jordan *-homomorphism
multiplied
by a fixed unitary element of $\li$. Consequently, we can assume that
our approximately local surjective isometry $\Phi$ is a unital Jordan
*-homomorphism. Now, the proof can be completed similarly to the
case of automorphisms, using \cite[Theorem 3]{MolStud1} in the place of
\cite[Theorem 2]{MolStud1}.
\end{proof}

To conclude the paper, we note that we feel it to be an interesting
algebraic-topological question to investigate for which
completely regular spaces $X$ holds it true that the operator algebra
$\BH$ "cleans up" the Stone-\v Cech compactification of $X$ in the
sense that we have seen above.


\end{document}